\documentclass[11pt,a4paper]{article}

\usepackage{amsmath}
\usepackage{amsfonts}
\usepackage{amssymb}
\usepackage{float}
\usepackage{float}
\usepackage{multicol}
\usepackage{graphicx}
\usepackage{subfig}
\usepackage{cases}
\usepackage{verbdef}
\usepackage{caption}
\usepackage[ruled,vlined]{algorithm2e}
\usepackage[font=small,labelfont=bf]{caption}

\title{Finding Geodesics on Surfaces Using Taylor Expansion of Exponential Map}
\author{E. Peyghan, E. Sharahi and A. Baghban}

\begin{document}
%%%%%%%%%%%%%%%%%%%%%%%%%%%%%
%%%%%%%%%%%%%%%%%%%%%%%%%%%%%
\maketitle
\begin{abstract}
Our aim in this paper, is to construct a numerical algorithm using Taylor expansion of exponential map to find geodesic joining two points on a 2-dimensional surface for which a Riemannian metric is defined. 
\end{abstract}
\section{Introduction}
After fixing the shape of manifold by a connection, the first steps are extracting the geodesics SODE and computing the curvature components. There are a wide variety of problems involving with geodesics in differential and computational geometry (e.g. \cite{berger, open, GABAY, Scheffer, udriste}). Actually, not only the geodesics are vastly applicable in various aspects of  science, but also they have a broad range. Among tons of works, we refer to \cite{farup, V1} in order to go through the heart of the matter and dealing with geodesics. But due to the fact that the nonlinear SODE which describes the geodesics has some analytical shortcomings in solving it, computational geometry, discretization and numerical analysis will appear (e.g. some basic effort in \cite{bose, fiori, ieee}). Geodesics of surfaces which start at a point in any direction, can be drawn by applying the numerical methods on the relative SODE for geodesics (e.g. \cite{peyghan}).  But, 
in what ways it is possible to find (maybe more than one) geodesics, which is about to join two arbitrary points on the surface (albeit if there exist)? Due to the by products of this question, any proper answer is welcomed. Here, to find  geodesics in between the points a method is described which is one of the many beneficial advantages of combining differential geometry and numerical analysis to come across geodesics in between the points. Implementing this algorithm has some profits. Indeed, it is simple, briefly sketched and it can be used in every problem in dealing with navigation, short paths, cost functions, manifold learning and maybe other stimulated issues. The main idea is the exponential map. Indeed, drawing the geodesic starting from $p\in M$ and ending to $q$ needs to know geodesic initial direction at $p$. Accordingly,  one can use the Taylor expansion of $\gamma(t)=\mathrm{exp}(tv)$ in local coordinate system
 \begin{align*}
 \gamma(t)=(p_1+v_1t+a_1^2t^2+a_1^3t^3+\cdots \\
 ,  p_2+v_2t+a_2^2t^2+a_2^3t^3+\cdots),
 \end{align*}
 where, using the geodesic equation 
 \begin{align*}
(\gamma^i)''(t)+\Gamma^i_{jk}(\gamma^j)'(t)(\gamma^k)'(t)=0,
\end{align*}
the coefficients $a^i_j$ are functions of $v^1,v^2$. So, if $\gamma(1)=q$ then we have a system of equations with respect to the $v^1$ and $v^2$:
\begin{numcases}{}
q^1=p_1+v_1t+a_1^2t^2+a_1^3t^3+\cdots, \nonumber\\
q^2=p_2+v_2t+a_2^2t^2+a_2^3t^3+\cdots, \nonumber
\end{numcases}
and solving the above system gives us the direction $v$. But, we have to use the numerical methods to solve. In this paper, an illustration of the algorithm is done. This algorithm is completely differential-geometric also it is easy to be implemented in problems dealing with $2$-dimensional manifolds imbedded in $\mathbb{R}^3$. 
Then after, obtaining geodesics in between two points in a Riemannian surface which is based on a Taylor expansion of exponential map, becomes significant. 
Back to our knowledge, it would happen to find more answers especially when the scalar curvature is positive. But one do not forget that in dealing with a special problem, requirement have to be prepared in various cases. So, it turns to modify the algorithms based on the advantages of the problem while there is no bewildering.
%
%It is remarkable that one can use all of these methods in the Finslerian cases only by a small change in geodesics equations. Also, for a complete discussion about geodesics in Finsler geometry, see \cite{shen1}.
%%%%%%%%%%%%%%%%%%%%%%%%%%%%%%%%%%%%
%%%%%%%%%%%%%%%%%%%%%%%%%%%%%%%%%%%%
\section{Preliminaries}
Through providing a manifold by an (affine) connection makes it possible to fix some of the features like curvature and minimal trajectories on that. It is well-known that any paracompact manifold, has a Riemannian metric. Let \((M^n, g)\) be a smooth Riemannian manifold, where $\nabla$ is devoted for its unique Levi-Civita connection. 
A vector field $X$ is parallely transported along a smooth curve $\gamma$, when $\nabla_{v_\gamma}X = 0$; where $v_\gamma$ is tangent to $\gamma$. If a smooth curve \(\gamma: [a, b] \subseteq \mathbb{R} \longrightarrow M\) parallel-transports its own tangent vectors, then it is a geodesic on \(M\). This means that \( \nabla_{\dot{\gamma}}\dot{\gamma} = 0\). It is easy to compute that the latter equality leads to the following system of ODE's
\begin{align}\label{main}
\dfrac{\mathrm{d}^2\gamma^k(t)}{\mathrm{d}t^2} + \Gamma^{k}{}_{ij}(t)\dfrac{\mathrm{d}\gamma^i(t)}{\mathrm{d}t}\dfrac{\mathrm{d}\gamma^j(t)}{\mathrm{d}t} = 0,
\end{align}
where \(  k \in \{1, \cdots, n\}\) and \(\gamma^i\)s are the components of \(\gamma\) and 
$\nabla_{\frac{\partial}{\partial x^i}}\frac{\partial}{\partial x^j} = \Gamma^{k}{}_{ij}\frac{\partial}{\partial x^k}$. The differential system of geodesic has the homogeniety property. Indeed, if $\gamma(t)$ is a geodesic, then for any nonzero constant $\lambda$, the curve $\gamma(\lambda t)$ is also a geodesic. Let $X_p \in \mathrm{T}_pM$ and suppose there exist a geodesic $\gamma_{X_p}: [0, 1] \rightarrow M$ satisfies  
\[
\gamma_{X_p}(0) = p, \ \ \frac{\mathrm{d}\gamma_{X_p}}{\mathrm{d}t}(0) = X_p. 
\]
Thus, the point $\gamma_{X_p}(1)$ is called the exponential of ${X_p}$ and denoted by $\mathrm{exp}_p({X_p})$.  Moreover, $\mathrm{exp}_p(t{X_p}) = \gamma_{t{X_p}}(1) = \gamma_{{X_p}}(t)$ and $\mathrm{exp}_p(0) = p$.
One of the most important consequences of the famous Hopf-Rinow theorem is that if \(M\) is a compact and connected Riemannian manifold, then any two points in \(M\) can join by a length minimizing geodesic. Discussions with lots of details, could be find in almost all Riemannian geometry books. For example, one can see a complete survey in \cite{postnikov}. 
Now, upon the intuition benefits of \(\mathbb{R}^3\),  we pay our attention to the case of surface \(S\) included in \(\mathbb{R}^3\). Consider such a surface with arbitrary coordinate system $(x,y)$ and with an arbitrary  Riemannian metric 
\begin{align*}
g=E\mathrm{d}x \otimes \mathrm{d}x +F\mathrm{d}x \otimes \mathrm{d}y +F \mathrm{d}y\otimes \mathrm{d}x + G \mathrm{d}y \otimes \mathrm{d}y,
\end{align*}%\label{0}
where \(E, F, G\in C^{\infty}(S)\). Using the well-known formula for Christoffel symbols 
$$\Gamma _{ij}^k=\frac{1}{2}g^{kl}\{\frac{\partial g_{li}}{\partial x ^j}+\frac{\partial g_{lj}}{\partial x ^i}-\frac{\partial g_{ji}}{\partial x ^l}\},$$
the system \eqref{main} translates to 
\begin{footnotesize}
\begin{numcases}{}
\begin{split}
&\alpha '' + \frac{GE_x-2FF_x+FE_y}{2(EG-F^2)}(\alpha ')^2+\frac{GE_y-FG_x}{EG-F^2}\alpha ' \beta '\\
&+\frac{2GF_y-GG_x-FG_y}{2(EG-F^2)}(\beta ' )^2=0,
\end{split} \nonumber \\
\begin{split}
&\beta '' +\frac{2EF_x-EE_y-FE_x}{2(EG-F^2)}(\alpha ')^2+\frac{EG_x-FE_y}{EG-F^2}\alpha ' \beta'\\
&+\frac{EG_y-2FF_y+FG_x}{2(EG-F^2)}(\beta')^2=0,
\end{split} \nonumber
\end{numcases}
\end{footnotesize}
where we supposed that \(\gamma = (\alpha, \beta)\) with respect to the coordinate system $(x,y)$.
\section{Constructing the algorithm}
We begin this section by assuming that $\gamma(t)=(x^1(t), x^2(t))$ is a geodesic passing through $p=(p^1,p^2)$ at $t=0$ with initial velocity vector $A=(a^1,a^2)$. Then by Taylor expansion, it follows 
\begin{align}\label{ixai}
x^i(t) = p^i + a^it + \sum_{n = 2}^{\infty}c^{i}_{n}t^{n},
\end{align}
for $i=1,2$.  Note that in this section we suppose that $(x^1,x^2)$ is our coordinate system. So, by calculating $c^i_j$s, we have the exponential mapping. From \eqref{ixai}, we find 
\begin{align}\label{final}
c^i_r = \frac{1}{r!} \stackrel{(r)}{{x}^i} (0) , \quad r = 2, 3,  \cdots.
\end{align}
Thus, to find $c^i_r$ generally, we should differentiate $r - 2$ times from \eqref{main}. Since $\Gamma ^i_{jk}(\gamma(t))=\Gamma ^i_{jk}(x^1(t), x^2(t))$, so
\begin{align*}
\frac{\mathrm{d}}{\mathrm{d}t}\Gamma ^i_{jk}(\gamma(t)) = \sum _{r=1}^{2} \frac{\partial \Gamma ^i_{jk}}{\partial x^r} (\gamma(t)) \dot{x}^r(t).
\end{align*}
Now, suppose that $q=(q^1,q^2)$ is a point in $\mathbb{R}^2$ and close to $p$. We know that the exponential mapping $\mathrm{exp}:\mathrm{T}_pS\to \mathbb{R}^2$ is locally diffeomorphism. So, we suppose that there exists a vector such as $A=(a^1,a^2)\in \mathrm{T}_pS$ so that $\mathrm{exp}(A)= q$. Now it is time to use the numerical methods to find such a vector $A$. If we find this vector then we can approximately find  the geodesic joining $p$ to $q$. This geodesic is $\gamma(t)=\mathrm{exp}(tA)$. So, we use Taylor expansion of $\gamma(t)=\mathrm{exp}(tA)$ at $t=1$ and solve the following system with respect to $A$
\begin{align}\label{ga}
\left\lbrace
\begin{array}{c}
p^1+a^1 +c^1_2+ \cdots + c^1_{n_0} = q^1 \\ 
p^2+a^2 +c^2_2+  \cdots + c^2_{n_0} = q^2,
\end{array} 
\right.
\end{align}
where $n_0+2$ is the order of our expansion. By solving the above system  we will get a geodesic starting from $p$ with endpoint   very close to point  $q$.
%%%%%%%%%%%%%%%%%%%%%%%%%%%%%
\subsection{Algorithm}
The Algorithm \emph{1} describes how we can do the above discussions by computer. Indeed, knowing values of $x(0)$ and $y(0)$ (as $x^1(0)$ and $x^2(0)$) and assuming $x^{(1)}(0)$ and $y^{(1)}(0)$ are independent variables, we use following procedure to approximate the geodesic using a polynomial of order of $n$. After finding all of the $n$ roots, the answers by minimum norm are geodesics.  
\begin{algorithm}
%	\captionsetup[algorithm]{name=Adams-Moulton algorithm}
%	\SetAlgoLined
%	\KwResult{Write here the result }
%	initialization\;
\For{$i = 2 : n$}{
Evaluate $x^{(i)}(0)$ and $y^{(i)}(0)$ using $x^{(i - 1)}(0)$, $y^{(i - 1)}(0)$, $x^{(i - 2)}(0)$ and $y^{(i - 2)}(0)$. \\
	Evaluate $c^1_i$ and $c^2_i$ using Eq.\eqref{final}.
}
Calculate $x^{(1)}(0)$ and $y^{(1)}(0)$ using a set of nonlinear equations as Eq.\eqref{ga}.
\For{$i = 2 : n$}{
	Evaluate numerical value of $x^{(i)}(0)$ and $y^{(i)}(0)$ using $x^{(i - 1)}(0)$ and $y^{(i - 1)}(0)$. \\
	Evaluate numerical value of $c^1_i$ and $c^2_i$ using Eq.\eqref{final}.
}
Calculate the numerical approximation of geodesic using Eq.\eqref{ixai}.
\caption{Finding Taylor expansion coefficients for the exponential map}
\end{algorithm}
Here, there are some famous examples that we examine the method on them as some practical implementations to ensure the accuracy of the algorithm. All of the three examples are classical and so, one can see that the algorithm find the geodesics with a high exactness.
\begin{figure}[H]%
    \centering
    \subfloat{{\includegraphics[width=5cm]{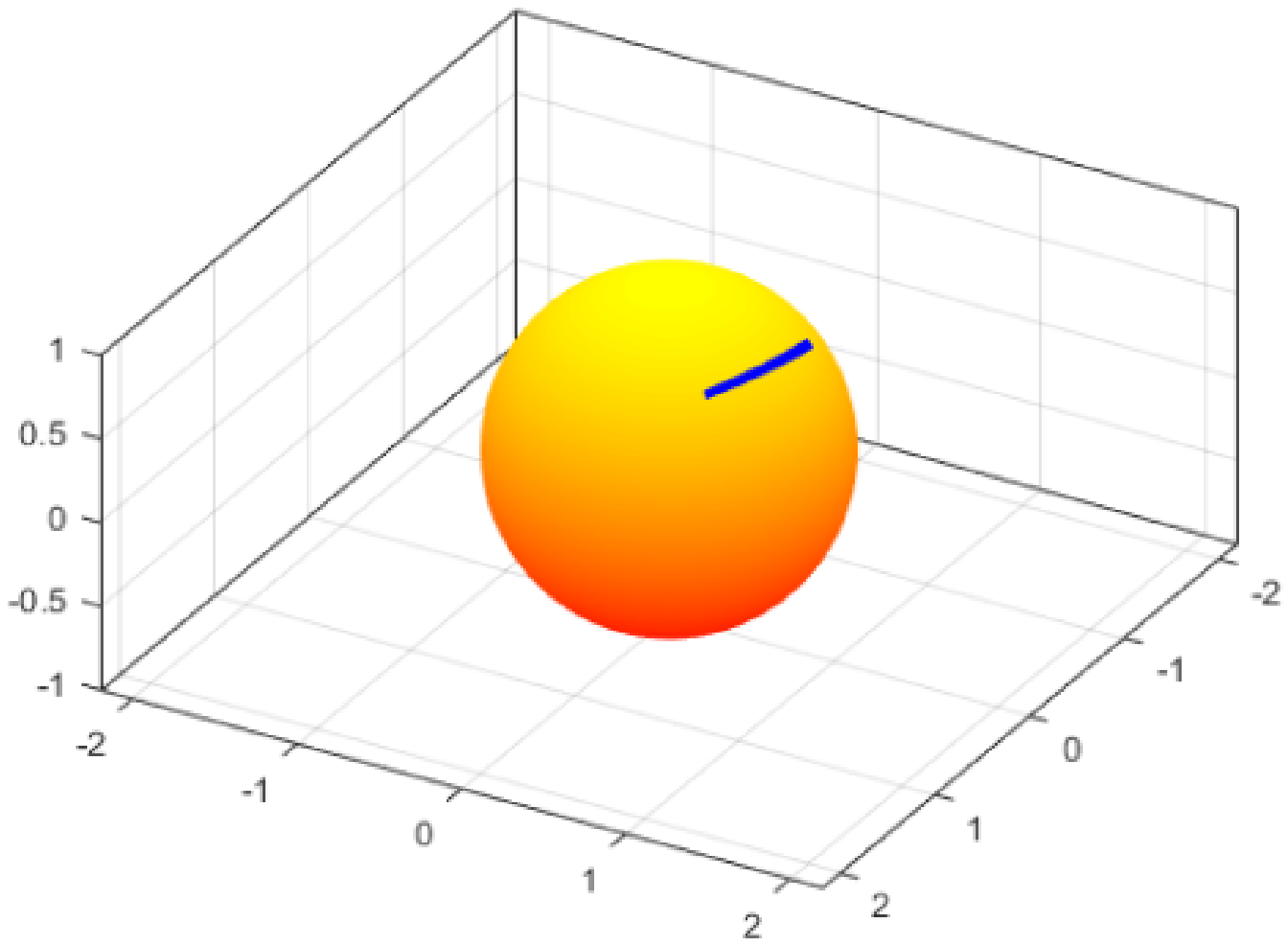} }}%
    \qquad
    \subfloat{{\includegraphics[width=5cm]{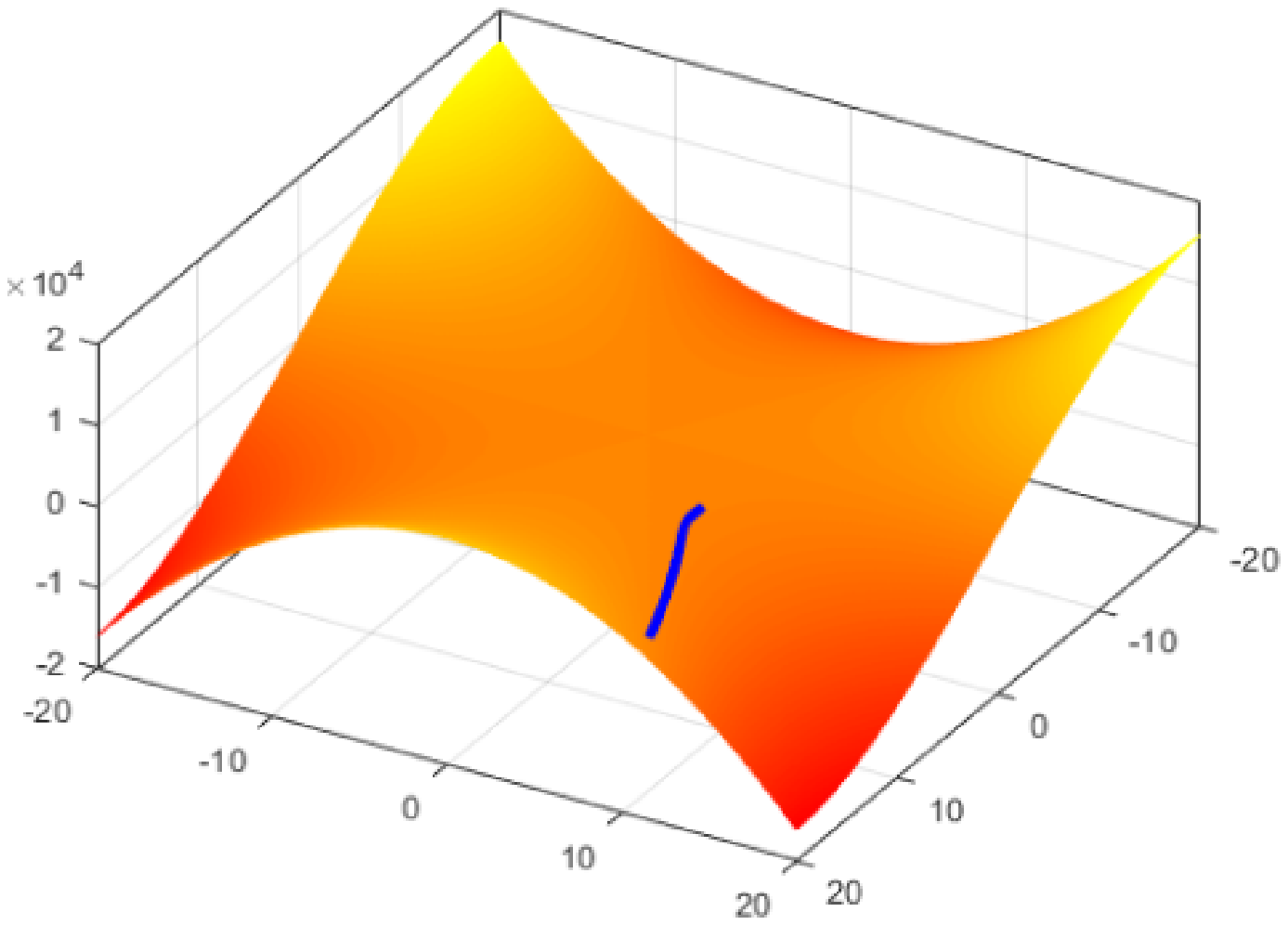} }}%
\end{figure}
\begin{tiny}
\begin{center}
\begin{tabular}{|c|c|c|c|c|}
\hline 
Surface & Initial point $p$ & End point $q$ & End point by the method & Order of expansion \\ 
\hline 
$z = \sqrt{1-x^2 - y^2}$ & $(1/2,1/2)$ & $(-1/3,2/3)$ & (-0.333333333  , 0.666666666) & 7 \\ 
\hline 
$z = x^3 - 3 x y^2$ & $(1,2)$ & $(15,7)$ & (14.999999987,  6.999997216) & 7 \\ 
\hline 
$\{ (x, y) | y > 0 \}$ & $(0.5,0.5)$ & $(0.55,0.6)$ & (0.549999999,  0.599999999) & 7 \\
\hline
\end{tabular} 
\end{center}
\end{tiny}
%
%%%%%%%%%%%%%%%%%%%%%%%%%%%

%%%%%%%%%%%%%%%%%%%%%%%%%%%%% 

Department of Mathematics, Faculty  of Science,\\ 
Arak University, Arak 38156-8-8349,  Iran. \\
\verb|epeyghan@gmail.com, esasharahi@gmail.com|

Department of Mathematics, Faculty  of Science, \\
Azarbaijan Shahid Madani University, Tabriz 53751 71379, Iran. \\
\verb|amirbaghban87@gmail.com|
\end{document}